\documentclass[12pt]{amsart}
\usepackage{amssymb,latexsym}
\newcommand{\vs}{\vspace{2mm}}
\theoremstyle{plain}
\newtheorem{theorem}{Theorem}
\newtheorem{remark}{Remark}
\newtheorem{corollary}{Corollary}

\newtheorem{lemma}{Lemma}
\newtheorem{proposition}{Proposition}
\theoremstyle{definition}

\newtheorem{claim}{Claim}
\newtheorem{example}{Example}
\theoremstyle{remark}

\numberwithin{equation}{section}

\begin{document}

\title[some galois extensions of quadratic extensions]{Some Galois extensions of quadratic extensions associated with Witt rings}
\author{Wenfeng Gao}
\address{537 Bellevue Way SE \\ Apt. 110 \\ Bellevue, WA 98004}
\email{wgao@sprintmail.com}
\author{J\'{a}n Min\'{a}\v{c}}
\address{The University of Western Ontario, Department of Mathematics \\ Middlesex College, London, Ontario, Canada N6A 5B7}
\email{minac@uwo.ca}
\thanks{The second author's research is partially supported by NSERC. This author also gratefully acknowledges the special
Dean of Science fund at The University of Western Ontario, along with the hospitality of the Mathematical Sciences Research
Institute in Berkeley, California during the Fall of 1999.} 
\date{}
\begin{abstract}
A Galois field extension $E/F$ whose Galois group is a pro-$2$-group of an exponent of at most $8$, with a nilpotency class
of at most $4$, is determined, such that it contains essential information about the Witt rings of all quadratic extensions
of $F$.
\end{abstract}

\maketitle

Let $F$ be a field of characteristic not $2$ and $F_q$ the quadratic
closure of $F$.  Let $G_F$ be the Galois group of $F_q/F$.  Then $G_F$
is a pro-$2$-group.  Let $F = F^{(1)}\subset F^{(2)}\subset\cdots\subset
F_q$ be a tower of fields such that $F^{(n+1)}$ is the compositum of
all of the quadratic extensions of $F^{(n)}$ which are Galois over $F$.
Observe that $F^{(n)}$ is a Galois extension of $F$.  We denote
$Gal(F^{(n)}/F)$ by $G^{[n]}_F$.  It is a quotient of $G_F$ and we
denote $G^{(n)}_F$ the kernel of the natural projection of $G_F$ onto $G_{F}^{[n]}$. We have 
$G^{(n+1)}_F=(G^{(n)}_F)^2[G^{(n)}_F,G_F]$. This means that $G^{(n+1)}_F$ is the
smallest closed subgroup of $G_{F}$ which contains all squares of elements in $G_{F}^{(n)}$ and all commutators $[\sigma,\tau]$, 
where $\sigma \in G_{F}^{(n)}$ and $\tau \in G_{F}$. We denote by $WF$ the Witt ring of a quadratic form over $F$. For the 
basic theory of Galois pro-$2$-extensions see [Ko], [N-S-W], [Rib], or [Sh]. For basic information on the Witt rings of quadratic
forms see [La] or [Sch]. We also use Kummer theory which describes abelian field extensions of a given exponent with enough
roots of unity in the base field. (See e.g. [A-T], Chapter 6.)

\vspace{2mm}

We use freely some well-known and easily derived conditions for the existence of an embedding of quadratic extensions in
Galois extensions with a Galois group isomorphic to a cyclic group of order $4$, and also conditions for an embedding of
biquadratic extensions in Galois extensions with a Galois group isomorphic to a dihedral group of order $8$, which we denote
as $D_{4}$. (See e.g. [G-M], pages 101 and 102.)

\vs

In [M-Sp] it was proved:

\begin{theorem}[M-Sp]                           
Let $F,F^\prime$ be two fields. Then

\noindent\begin{enumerate}
\item[(1)] {$WF\cong WF^\prime$ implies that $G^{[3]}_F\cong G^{[3]}_{F^\prime}$,}
\item[(2)] {If $\langle 1,1\rangle_F$ is universal, assume that $\sqrt{-1}\in F$ iff $\sqrt{-1}\in F^\prime$. Then   
$G^{[3]}_F\cong G^{[3]}_{F^\prime}$ implies that $WF\cong WF^\prime$.}
\end{enumerate}

\end{theorem}

This theorem shows that the theory of Witt rings can be thought of as part
of Galois theory.  In particular the classification of Witt rings can
be reduced to the classification of Galois groups $G^{[3]}_F$ and specifying whether
$\sqrt{-1}\in F$ or $\sqrt{-1}\notin F$.

\vspace{2mm}

One of the major unsolved probems in the theory of quadratic forms is
the precise relationship between $WF$ and $WF(\sqrt{a})$, where 

$$a\in F^* \colon = F -\{0\} {\mbox{ and }} F(\sqrt{a}) {\mbox{ is a quadratic extension of}} \; F.$$

However some interesting progress has been done on some special cases,
and there is partial information related to progress with a general case. (See e.g., [E-L], [L-Sm], [M-Sm], [M-W],
[P-S-C-L] and [Szy].) Closely related to this problem are the cohomology rings of $G_{F}^{[3]}$ which can be viewed as
invariants of Witt rings. (See [A-K-M].)

\vspace{2mm}

We set $L = F(\sqrt{a})$ such that $a\in F^{*}$ and $[L:F] = 2$.  Then we
have $F\subset L\subset F^{(2)}$. Also for any pro-$2$-extension $K/F$ we denote by $\bar{K}/F$ or simply by
$\bar{K}$ the Galois closure of $K/F$. In our paper we clarify ``how much Galois theory'' is needed to understand $WL$.
More precisely we prove

\vspace{3mm}

\begin{theorem}
$\overline{L^{(3)}}\subset F^{(5)}$.
\end{theorem}

Therefore the Galois closure of the compositum of all extensions $L^{(3)}$ is a subfield of $F^{(5)}$. We also show that in
general the Galois closure of the compositum of all extensions $L^{(3)}$ is not a subfield of $F^{(4)}$. (See Example 1 
below.) Hence we see that this $F^{(5)}$ in Theorem 2 above cannot be replaced by $F^{(4)}$. Finally we shall observe in 
Proposition 1 below that $L^{(3)}$ is Galois over $F$ and therefore the bar over $L^{(3)}$ in \linebreak Theorem 2 may be omitted. We 
conclude our paper with the observation that the Galois group of the compositum of all $L^{(3)}$ where $L$ runs over 
quadratic extensions of $F$, has an exponent of at most $8$. Therefore this compositum is in general not equal to $F^{(5)}$. 
(See Proposition 2.)

\vspace{2mm}

We denote $H^1(G^{(2)}_F,\mathbb{F}_{2})^{G^{[2]}_F}$ by $J_1(G^{(2)}_F)$. Thus $J_1(G^{(2)}_F)$ denotes the fixed elements
of $H^{1}(G_{F}^{(2)},\mathbb{F}_{2})$ under the natural action of $G_{F}^{[2]}$. Here $H^{1}(G_{F}^{(2)}, \mathbb{F}_{2})$
is the first cohomology group of $G_{F}^{(2)}$ with coefficients in $\mathbb{F}_{2}$. Thus by Kummer theory $J_{1} (G_{F}^{(2)}) 
\cong ((F^{(2)})^{*} / ((F^{(2)})^{*})^{2})^{G_{F}^{[2]}}$, where again the latter group denotes the fixed elements of the
square class group of the multiplicative group of $F^{(2)}$. (See [A-T], Chapter 6, \linebreak Section 2 for details on Kummer theory,
and [G-M], page 100, for a more detailed justification of the isomorphism above. In the rest of the paper we identify the
two isomorphic groups mentioned above.) Observe that $F^{(3)} \colon = F^{(2)}(\sqrt{J_{1}(G_{F}^{(2)}}))$ where the latter
field means a compositum of all fields $F^{(2)}(\sqrt{\gamma}),\gamma \in J_{1}(G_{F}^{(2)})$. Similarly we have
$F^{(n+1)} = F^{(n)}(\sqrt{J_1(G_{F}^{(n)}}))$ where $J_{1}(G_{F}^{(n)})$ denotes the fixed elements of 
$H^{1}(G_{F}^{(n)},\mathbb{F}_{2})$ under the action of $G_{F}^{[n]}$ for each $n \geq 2$. 

\vspace{2mm}

We also frequently use the following well-known fact. (See e.g. [War].)

\vspace{2mm}

\noindent{\bf Statement.} {\em Suppose that $K/F$ is a Galois extension, $G = Gal(K/F)$ is its Galois group
and $a \in K^{*}$. Then $K(\sqrt{a})/F$ is Galois iff $\sigma(a)/a$ is a square in $K^{*}$ for each $\sigma \in G$.} 

\vspace{2mm}

We begin with recalling some of the basic inclusions between the fields that we are dealing with. We shall use the notation
$\sigma/K$ to mean a restriction of an automorphism $\sigma$ of an overfield of $K$ to $K$. We also use the symbol $[b]_{K}$
to denote a class of $K^{*}/K^{*2}$ represented by $b$. Sometimes when there is no danger of confusion we simply write 
$[b]$ to denote $[b]_{K}$.   

\vspace{3mm}

\begin{lemma}
$L^{(2)}\subset F^{(3)}$.
\end{lemma}

\begin{proof}
By the definition of $L^{(2)}$, it is enough to show that for each quadratic extension
$L(\sqrt{b})/L$, $b\in L^*$, we have $L(\sqrt{b})\subset F^{(3)}$.
Consider any element $\sigma\in G^{[2]}_F$, we have $\sigma |_L\in
Gal(L/F) = \{1,\tau\}\cong\mathbb{Z}/2$.  If $\sigma |_L = 1$, then
$\sigma (b) = b$ and $\sigma (b)/b = 1\in (F^{(2)^*})^2$.  If $\sigma
|_L = \tau$, then $\sigma (b)/b = \sigma (b)b/b^2 = N_{L/F}(b)/b^2\in
(F^{(2)^*})^2$.  Therefore $\sigma (b)/b\in (F^{(2)^*})^2$ in each case
and consequently $[b]\in J_1(G^{(2)}_F)$.  This means that
$F^{(2)}(\sqrt{b})\subset F^{(3)}$.  Because $L(\sqrt{b})\subset
F^{(2)}(\sqrt{b})\subset F^{(3)}$, we have our desired inclusion.
\end{proof}

\begin{lemma}
$F^{(3)}\subset L^{(3)}$.
\end{lemma}

\begin{proof}
By the definition of $F^{(3)}$, it is enough to show that for each $[b]\in
J_1(G^{(2)}_F)$ we have $F^{(2)}(\sqrt{b})\subset L^{(3)}$.  Consider
any $\sigma\in Gal(L^{(2)}/L)$.  We have $\sigma (b)/b =
\sigma|_{F^{(2)}}(b)/b\in (F^{(2)^*})^2\subset (L^{(2)^*})^2$, the first
inclusion follows from the inclusion $[b]\in J_1(G^{(2)}_F)$.  Therefore $[b]_{L^{(2)}}\in
J_1(G^{(2)}_L)$ as well.  Hence $F^{(2)}(\sqrt{b})\subset
L^{(2)}(\sqrt{b})\subset L^{(3)}$ for each $[b]\in J_1(G^{(2)}_F)$ which
proves that $F^{(3)}\subset L^{(3)}$ as we claimed.
\end{proof}

Therefore we have the tower $F\subset L\subset F^{(2)}\subset
L^{(2)}\subset F^{(3)}\subset L^{(3)}$.

\vspace{2mm}

Let $K/F$ be a Galois 2-extension such that $L\subset K$. This means that the Galois group of $K/F$ is a pro-$2$-group. 
Let $T/K$ be a quadratic extension such that $T/L$ is Galois.  Let $\overline{T}$ be the Galois closure of $T$ over $F$.

\begin{lemma}
There are only two possibilities:
\begin{enumerate}
\item[(1)]$T$ is Galois over $F$.  Then $T = \overline{T}$.
\item[(2)]$T$ is not Galois over $F$.  Then $\overline{T} = T(\sqrt{t})$ for some $t\in T^{*}- \, T^{*2}$.
\end{enumerate}
\end{lemma}

\begin{proof}
It is enough to show that if $T$ is not Galois over $F$ then
$\overline{T}/F$ is a quadratic extension of $T$.  Assume that $T/F$ is not Galois.  Set $T = K(\sqrt{k})$, $G = Gal(K/F)$
and $H = Gal(K/L)$. Because $T/F$ is not Galois we have $\sigma(k)/k\notin (K^*)^2$ for some $\sigma \in G-H$. We set 
$E = T(\sqrt{\sigma(k)})$. We want to show that $E$ is Galois over $F$.  Consider any $\tau\in G_{F}$. If we can show that 
$\tau(E)\subset E$, we shall be done.  Because $K/F$ is Galois and $E = K[\sqrt{k},\;\sqrt{\sigma(k)}]$ as a ring,
we see that it is enough to show that $\tau(\sqrt{k})),\;\tau(\sqrt{\sigma(k)}) \in E$.

\vspace{2mm}

Consider first $\tau(\sqrt{k})$.  We have $\tau (k) = \tau ((\sqrt{k})^2) = (\tau(\sqrt{k}))^2$.  Therefore 
$\tau (\sqrt{k})$ must be some square root of $\tau (k)$.  If $\tau |_L = 1$, then $\tau |_K\in H$ and 
$\tau (k)/k\in (K^*)^2$.  Therefore $\sqrt{\tau (k)}\in T$. Suppose now that $\tau |_L\not= 1$.  Since 
$(\tau^{-1}\sigma)|_{L} = 1$, then $[\tau^{-1} \sigma (k)]_K = [k]_K$.  Therefore $[\tau (k)]_K = [\sigma (k)]_K$.  But this
means that $\sqrt{\tau (k)}\in E$.

\vspace{2mm}

In a similar fashion, by exploiting the fact that the orbit of the element $[k]$ of $K^*/(K^*)^2$ under the action of $G$ 
has $2$ elements we prove that $\tau(\sqrt{\sigma(k)}) \in E$ for each element $\tau \in G$. Therefore $E/F$ is Galois.  
Since $T/F$ is not Galois and since $[E:T] = 2$. Therefore we see $E = \bar{T}$.  Our proof is complete.
\end{proof}

\begin{lemma}
$(F^{(4)}\cap L^{(3)})/F$ is Galois.
\end{lemma}

\begin{proof}
Both extensions $F^{(4)}$ and $L^{(3)}$ are Kummer extensions of $F^{(3)}$ of exponent at most $2$. From our inductive
definition of $F^{(4)}$ using $F^{(3)}$ we see that each subextension between $F^{(3)}$ and $F^{(4)}$ is Galois over $F$.
In particular, $L^{(3)}\cap F^{(4)}$ is Galois over $F$.
\end{proof}

\begin{proof}[Proof of Theorem 2.]
By the definition of $L^{(2)}$ and Lemma 1, it is enough to show that for each $\gamma\in L^{(2)^*}$
such that $L^{(2)}(\sqrt{\gamma})/L$ is Galois; hence we have a Galois closure of $(F^{(4)}\cap L^{(3)})(\sqrt{\gamma})$
over $F$, which is contained in $F^{(5)}$. (Observe that $L^{(2)} \subset F^{(4)} \cap L^{(3)}$ and therefore $L^{(2)}(\sqrt{\gamma})
\subset (F^{(4)} \cap L^{(3)})\sqrt{\gamma}.)$ If $\sqrt{\gamma}\in F^{(4)}\cap L^{(3)}$, then
$L^{(2)}(\sqrt{\gamma})\subset F^{(4)}\subset F^{(5)}$.  Otherwise, we can apply Lemma 3 to the following situation:
$K = F^{(4)}\cap L^{(3)}$ and $T = (F^{(4)}\cap L^{(3)})(\sqrt{\gamma})$. Indeed $K/F$ is Galois by Lemma 4, and since 
$L^{(2)}(\sqrt{\gamma})$ is Galois over $L$ and $L^{(2)} \subset F^{(4)} \cap L^{(3)}$, we see that $T$ is also Galois over
$L$. We also use notation employed in the proof of Lemma 3; in particular $\sigma \in G := Gal(K/F)$ and       
$\sigma /L \neq 1_{L}$. For the sake of ease of notation we denote by $\sigma$ also its restriction to $F^{(3)}$. 
Then we obtain $\overline{T} = T(\sqrt{\gamma}, \sqrt{\sigma(\gamma)})$. Moreover we have $\gamma\in L^{(2)}\subset F^{(3)}$,
so $\gamma\sigma(\gamma)\in F^{(3)}$.  We have $\tau (\gamma)\gamma^{-1}\in (L^{(2)^*})^2$ for each $\tau$ which fixes $L$.
Therefore for each $\tau \in Gal(K/L)$ restricted to $F^{(3)}$, we have $\tau\sigma = \sigma\tau^\prime$, $\tau^\prime\in G^{[3]}_F$ and $\tau^\prime |_L = 1$
and 
\[ \tau (\gamma\sigma (\gamma )) = \tau (\gamma )\sigma (\tau^\prime(\gamma ))\equiv\gamma\sigma (\gamma)\;\quad
\mbox{mod}\; (F^{(3)^*})^2\; . 
\]
On the other hand, if $\delta = \sigma\tau^\prime$, where $\tau^\prime\in Gal(F^{(3)}/L)$ we see:
\[ \delta (\gamma\sigma (\gamma )) = \sigma (\tau^\prime (\gamma
))\sigma\tau^\prime\sigma (\gamma )\equiv\sigma (\gamma
)\gamma\quad\mbox{mod}\; (F^{(3)^*})^2\; . \]
Therefore we see that $F^{(3)}(\sqrt{\gamma\sigma(\gamma)}) \subset 
F^{(4)}$.  But this means that $\gamma\equiv\sigma(\gamma)$ mod
$(F^{(4)^*})^2$ and therefore $\gamma\in J_1(G^{(4)}_F)$.  Thus we
obtain $\overline{L^{(3)}}\subset F^{(5)}$ as we claimed.
\end{proof}

\begin{remark}
{\textnormal{Observe that in the proof above we also showed the following statement:}}

\vs

Let $\gamma$ be an element of $L^{(2)*}$ such that $L^{(2)}(\sqrt{\gamma})/L$ is a Galois extension. Let also 
$\sigma\in G_{F}^{[3]}-Gal(F^{(3)}/L)$. Then the Galois closure of the field extension $L^{(3)}\cap F^{(4)}(\sqrt{\gamma})/F$ 
is a field $T = L^{(3)} \cap F^{(4)}(\sqrt{\gamma}, \sqrt{\sigma(\gamma)})$ and there is a tower of field extensions which
are Galois over $F$: $L^{(3)} \cap F^{(4)} \subset L^{(3)} \cap F^{(4)}(\sqrt{\gamma\sigma(\gamma)}) \subset T$. Moreover
$L^{(3)}\cap F^{(4)}(\sqrt{\gamma\sigma(\gamma)}) \subset F^{(4)}$ and $T \subset F^{(5)}$. 
\end{remark}

\begin{corollary}
The compositum of all $W$-fields of all quadratic extensions of $F$ is
in $F^{(5)}$.
\end{corollary}

Observe also the following simple corollary of our discussion. 

\begin{proposition}
Suppose that $L/F$ is any quadratic extension of $F$. Then $L^{(3)}/F$ is a Galois extension. 
\end{proposition}

\begin{proof}
Suppose that $\tau \in G_{F}$. It is sufficient to prove that $\tau(L^{(3)}) \subset L^{(3)}$. First observe that 
$\tau(L^{(2)}) \subset L^{(2)}$ as $\sqrt{\tau(l)} \in L^{(2)}$ for each $l \in L^{*}$. 

Now let $\gamma \in L^{(2)*}$ be an element of $L^{(2)*}$ such that $L^{(2)}(\sqrt{\gamma})/L$ is Galois. Then for each 
$\delta \in Gal(L^{(2)}/L)$ we have $\delta(\gamma)\gamma \in (L^{(2)})^{*2}$. Now let us consider two cases. In the second
case we use the fact that $G_L$ is the normal subgroup of $G_F$ of index $2$. 

\vs

{\bf Case 1} $\tau/L$ is an identity automorphism. Then for each 

$$\omega \in Gal(L^{(2)}/L) {\mbox{ we have }} [\omega \tau(\gamma)]=[\tau(\gamma)]=[\gamma] {\mbox{ in }} (L^{(2)})^{*}/(L^{(2)*})^{2}.$$ 

\vs

{\bf Case \,2} $\,\tau/L$ is a nontrivial automorphism. Then for each $\omega\in Gal(L^{(2)}/L)$ we can find $\omega_{1}, 
\omega_2 \in G_L$ such that $\omega_1 \mid L^{(2)} = \omega$ and $\tau \omega_2 = \omega_{1}\tau$. Therefore 
$[\omega\tau(\gamma)]=[\omega_1\tau(\gamma)]=[\tau\omega_2(\gamma)]=[\tau(\gamma)]$ in $L^{2*}/(L^{(2)*})^2$. 

\vs

Hence in both cases we see that $\tau(\gamma)\in J_1(G_{L}^{(2)})$ and therefore $L^{(2)}(\sqrt{\tau(\gamma)})/L$ is a 
Galois extension. From the fact that $L^{(3)}$ is the compositum of all fields $L^{(2)}(\sqrt{\gamma}), \gamma \in 
L^{(2)*}$ such that $L^{(2)}(\sqrt{\gamma})/L$ is Galois, we see that $\tau(L^{(3)}) \subset L^{(3)}$. Our proof is now 
completed. 
\end{proof}

The following is an example of $L^{(3)} \nsubseteq F^{(4)}$. Here we shall employ the usual notation: $\mathbb{R}$ - the 
field of all real numbers - and $\mathbb{C}$ - the field of all complex numbers. 

\vs

\begin{example}
Let $F = \mathbb{R}(t)$. Then $L \colon = F(\sqrt{-1}) = \mathbb{C}(t)$. Since $\mathbb{C}[t]$
is a unique factorization domain, then $[t + 2i] \neq [1], [t-i] \neq [1] \in L^{*}/L^{{*}2}$. Also $[t+2i]$ and $[t-i]$ are
linearly independent in $L^{*}/L^{{*}2}$. It is well-known that $H^{2}(\mathbb{C}(t),\mathbb{F}_{2}) = 0$. (See e.g. [Rib], page
70 and page 217.) Then there exists a Galois extension $K/L$ such that $Gal(K/L)\cong D_{4}$ and $\mathbb{Z}/4\mathbb{Z} \cong 
Gal(K/L(\sqrt{(t+2i)(t-i)}))$.                               

\vs

We denote this $D_{4}$-extension of $L$ by $K$. We also denote 

$$\mathbb{C}(\sqrt{t-i},\sqrt{t+2i}) {\mbox{ by }} E.$$

Observe that

$$F^{(2)} = \mathbb{C}(\sqrt{t-r},\sqrt{t^{2} + bt + c}, b^{2} - 4c < 0, b,c,r \in \mathbb{R}).$$
\end{example}

\begin{claim}
$E \cap F^{(2)} = \mathbb{C}(t)$.
\end{claim}

\begin{proof}
Set $V$ as the subspace of $L^{*}/L^{{*}2}$ generated by $[t-r],[t^{2}+bt+c],r,b,c \in \mathbb{R}$ such that $b^2-4c < 0$. 
In order to prove our claim, it is enough to show that $[t-i]$ and $[t+2i]$ generate a vector space $W$ of $L^{*}/L^{{*}2}$ 
such that $W \cap V = \{0\}$. Assume that some element; say $[(t-i)(t+2i)] \in W \cap V$. Then 

$$(t-i)(t+2i)P(t)^{2} = \prod_{r}(t+r) \prod_{b,c}(t^{2}+bt+c)Q(t)^{2}$$

for some $P(t), Q(t) \in \mathbb{C}[t]$ and some finite set of linear factors $t+r$ and a finite set of quadratic factors
$(t^{2}+bt+c)$ such that $b^2-4c < 0$. This implies that $(t-i)$ must divide some $t^{2}+bt+c$. However if $t-i \mid t^{2}+bt+c$, 
then $t^{2}+bt+c = (t-i)(t+i) = t^{2}+1$. Hence we see that $t+i$ appears on the right hand side of our equation above, in
odd power, but $t+i$ appears on the left hand side of our equation above, in even power - a contradiction. Hence $E \cap F^{(2)} = L$.
\end{proof}

\begin{claim}
$K \cap F^{(2)} = L$.
\end{claim}

\begin{proof}
Indeed $K \cap F^{(2)} \subset L^{ab} \colon =$ the maximal abelian subextension of $K$ over $L$. (Because $F^{(2)}/L$ is
an abelian extension.) Hence $K \cap F^{(2)} \subset L^{ab} = E$. By {\it{Claim 1}}, we obtain 
$K \cap F^{(2)} \subset E \cap F^{(2)} = \mathbb{C}(t) = L$. 

By Claim 2, $Gal(KF^{(2)}/F^{(2)}) = Gal(K/K \cap F^{(2)}) \cong D_{4}$. Observe however that $F^{(4)}/F^{(2)}$ is an 
abelian extension. Indeed one can check easily that $[G_{F}^{(2)},G_{F}^{(2)}] \subset G_{F}^{(4)}$. (See e.g. Chapter 7.15 in
[Ko].) Then $KF^{(2)} \nsubseteq F^{(4)}$. On the other hand, $KF^{(2)} \subset L^{(3)}$. Hence $L^{(3)} \nsubseteq F^{(4)}$. 
Our Example 1 is completed. 
\end{proof}

\vs

Set $(\sqrt/F)^{(3)} :=$ the compositum of all field extensions $L^{(3)}/F$, where $L$ is a quadratic extension of $F$,
and $Gal((\sqrt/F)^{(3)}/F) \colon = G_{\sqrt/F}^{[3]}$. From Proposition 1 we know that $(\sqrt/F)^{3}$ is a Galois
extension. We proved that $(\sqrt/F)^{(3)} \subset F^{(5)}$. Observe however that in general the exponent of
$G_{F}^{[5]} = Gal(F^{(5)}/F)$ is $16$ while the exponent of $Gal((\sqrt/F)^{(3)}/F)$ is in general just $8$. More precisely,
we have the proposition below. Recall that $F$ is quadratically closed if $F_{q}=F$ and $F$ is a Euclidean field if 
$F_{q}^{*2}$ is an ordering of $F_{q}$ and $F_{q}^{*}=F_{q}^{*2}\cup -F_{q}^{*2}$. In the proof of the proposition below we
consider pythagorean fields. Recall that a field $F$ is pythagorean if each sum of two squares is again a square. It is well
known that if $F$ is a pythagorean field then $F(\sqrt{-1})$ contains all of the roots of unity of an order with a power of
$2$. (See e.g. [Be2], page 83.) 

\begin{proposition}
\vs
\begin{enumerate}
\item{If $F$ is quadratically closed, then $(\sqrt/F)^{(3)} = F = F^{(n)}$ for each $n \in \mathbb{N}$.}
\item{If $F$ is a Euclidean field then $(\sqrt/F)^{(3)}=F^{(2)}=F^{(n)}$ for each $n \in \mathbb{N}, \, n \geq 2$.}
\item{If $F$ is neither a quadratically closed nor a Euclidean field, then $(\sqrt/F)^{(3)} \subset F^{(5)}$ and the exponent
of $G_{F}^{[5]}$ is $16$ while the exponent of $Gal((\sqrt/F)^{(3)}/F)$ is $8$. In particular $(\sqrt/F)^{(3)} \neq F^{(5)}$.}
\end{enumerate}
\end{proposition}

\begin{proof}
If $F$ is a Euclidean field then $F_{q}=F(\sqrt{-1})$. (See [Be1], Satz 3.) Therefore both cases 1 and 2 are clear. 
Now consider the case when $F$ is neither quadratically closed nor a Euclidean field. Then from the definition of $F^{(5)}$
we see that $G_{F}^{[5]}$ has an exponent of at most $16$. If $F$ is not a pythagorean field, then there exists an element 
$a \in F^{*}-F^{*2}$ which is a sum of two squares. Then $F(\sqrt{a})$ can be imbedded in a Galois extension $L/F$
such that $Gal(L/F)\cong\mathbb{Z}/4\mathbb{Z}$. Then it is well-known (see e.g. [Ku-Le], Theorem 2) that for each 
$n \in \mathbb{N}$ a Galois extension $K/F$ exists such that $Gal(K/F)$ is $\mathbb{Z}/2^{n}\mathbb{Z}$. In particular 
a Galois extension $K/F$ exists such that $Gal(K/F)=\mathbb{Z}/16\mathbb{Z}$. Then $K \subset F^{(5)}$ and any element
$\sigma \in G_{F}^{[5]}$ which restricts to a generator of $Gal(K/F)$ has an order of $16$. Hence the exponent of $G_{F}^{[5]}$
is $16$. 

On the other hand if $F$ is a pythagorean field which is neither quadratically closed nor a Euclidean field, we can pick an
element $a \in F^{*}$ such that the set $\{[a],[-a]\}$ is linearly independent over $\mathbb{F}_{2}$. In particular we see
that $[a] \neq [1]$ and $a \notin -4F^{*4}$. Therefore from the well-known theorem (see e.g. [Lan], Theorem 9.1) we see that
the polynomial $X^{16}-a$ is irreducible over $F$. Since $F$ is a pythagorean field we know that $F(\sqrt{-1})$ contains all
$2^{n}-th$ roots of unity, $n \in \mathbb{N}$. Therefore the splitting field of the polynomial $f(X)=X^{16}-a$ is 
$N=F(\sqrt{-1},\sqrt[16]{a})$ where $\sqrt[16]{a}$ means any element $b \in F_{q}^{*}$ such that $b^{16}=a$. 

\vspace{2mm}

Observe that from Kummer theory we know that $Gal(N/F(\sqrt{-1}))\cong\mathbb{Z}/16\mathbb{Z}$ and that we can choose as a
generator of $Gal(N/F(\sqrt{-1}))$ an element $\tau$ such that $\tau(\sqrt[16]{a})=\zeta_{16}\sqrt[16]{a}$ where $\zeta_{16}$
is a primitive $16$th root of unity. Observe that $F(\sqrt{-1},\sqrt[4]{a})\subset F^{(3)}$ and therefore $F(\sqrt{-1},
\sqrt[16]{a})\subset F^{(5)}$. Therefore we see that for each $\delta\in G_{F}^{[5]}$ such that $\delta(\sqrt{-1})=\sqrt{-1}$
and $\delta/N = \tau$ we have an order of $\delta$ that is $16$. Thus again we see that the exponent of $G_{F}^{[5]}$ is $16$. 

\vspace{2mm}

Now consider $L^{(3)}$, for any $L=F(\sqrt{a}), a \in F^{*} - F^{*2}$ and any $\sigma\in G_{F}^{[5]}$. Then $\sigma^{2}\in
Gal(F^{(5)}/L)$ and we can consider the restriction of $\sigma^{2}$ to $L^{(3)}/L$. Since each element of $Gal(L^{(3)}/L)$
has an order of at most $4$, we see that $\sigma^{8}/L^{(3)}$ is the identity automorphism of $L^{(3)}$. Because 
$(\sqrt/F)^{(3)}$ is the compositum of all extensions $L^{(3)}/L$ over $F$ we see that the exponent of $Gal((L/F)^{(3)}/F)$
is at most $8$. 

\vspace{2mm}

To show that the exponent of $Gal((\sqrt/F)^{3}/F)$ is exactly $8$, we shall proceed in the same fashion as in the beginning of 
our proof. 

\vspace{2mm}

Namely if $F$ is not a pythagorean field, we see again that a Galois extension $K/F$ exists such that $Gal(K/F)\cong\mathbb{Z}/8\mathbb{Z}$. 
Let $L$ be the unique quadratic subextension of $K/F$. Then $Gal(K/L)\cong\mathbb{Z}/4\mathbb{Z}$ and we see that 
$K \subset L^{(3)}$. Now let $\sigma$ be any element of $Gal((\sqrt/F)^{3}/F)$ such that its restriction to $K/F$ is a
generator of $Gal(K/F)$. Then the order of $\sigma$ is at least $8$, and we can conclude that the exponent of our Galois
group $Gal((\sqrt/F)^{3}/F)$ is $8$. 

\vspace{2mm}

If $F$ is a pythagorean field which is neither a quadratically closed nor a Euclidean field, we can again find an element 
$a \in F^{*}$ such that the polynomial $g(X)=X^{8}-a$ is irreducible over $F$. Then the field extension $N=F(\sqrt{-1},\sqrt[8]{a})$
is the splitting field of $g(X)$ and $Gal(N/F(\sqrt{-1}))\cong\mathbb{Z}/8\mathbb{Z}$. Also observe that $N$ is the Galois
closure of $L(\sqrt[4]{(\sqrt{a})})$ over $F$ where $L=F(\sqrt{a})$. Therefore $N \subset (\sqrt/F)^{(3)}$. Thus we can again
conclude that the exponent of $Gal((\sqrt/F)^{(3)}/F)$ is $8$.  
\end{proof}

\end{document}